\documentclass[11pt]{amsart}
\usepackage{amsmath}
\usepackage{amscd}
\usepackage{amsthm}
\usepackage{graphicx}
\usepackage{amssymb}
\usepackage{epstopdf}
\theoremstyle{plain}

\theoremstyle{definition}

\numberwithin{equation}{section}

\title{On the extension problem of pluricanonical forms \\
\rm Dedicated to Professor Friedrich Hirzebruch on his seventieth birthday}
\author{Yujiro Kawamata}
\usepackage{fancyhdr}
\begin{document}

\maketitle

%\address Department of Mathematical Sciences, University of Tokyo, 
%Komaba, Meguro, Tokyo, 153-8914, Japan \endaddress
%\email kawamata\@ms.u-tokyo.ac.jp\endemail

%\keywords plurigenus, canonical singularity, deformation, 
%multiplier ideal sheaf, vanishing theorem 
%\endkeywords
%\subjclass 14B07, 14E30, 14F17
%\endsubjclass

%\abstract
We review some recent development on the extension problem of 
pluricanonical forms from a divisor to the ambient space
in [Si], [K5] and [N3] with simplified proofs.

A section for a correction is added.
%\endabstract

%\endtopmatter

%%%%%%%%%%%%%%%%%%%%%%%%%%%%%
%%%%%%%%%%%%%%%%%%%%%%%%%%%%%
%%%%%%%%%%%%%%%%%%%%%%%%%%%%%
\section{Introduction}

The purpose of this paper is to review some recent development on the 
extension problem of pluricanonical forms from a divisor to the ambient space.
The main tools of the proofs are the multiplier ideal sheaves and the
vanishing theorems for them.

Let $X$ be a compact complex manifold.
The {\it $m$-genus} $P_m(X)$ of $X$ for a positive integer $m$ is defined by
$P_m(X) = \text{dim }H^0(X, mK_X)$.
The growth order of the plurigenera for large $m$ is called the {\it Kodaira 
dimension} $\kappa(X)$: we have $P_m(X) \sim m^{\kappa(X)}$ for any 
sufficiently large and divisible $m$.  
We have the following possibilities: $\kappa(X) = -\infty, 0, 1, \cdots$, or
$\text{dim }X$.
In particular, if $\kappa(X) = \text{dim }X$, 
then $X$ is said to be of {\it general type}.
It is important to note that these invariants are independent of
the birational models of $X$.

The plurigenera are fundamental discrete invariants
for the classification of algebraic varieties.
But they are by definition not topological invariants.
However, in order that such classification is reasonable, the following 
conjecture due to Iitaka should be true:

\vskip 1pc

{\bf Conjecture 1.1.}\footnote{This is proved by Y.-T. Siu.}
{\em Let $S$ be an algebraic variety, and let
$f: \mathcal X \to S$ be a smooth algebraic morphism.
Then the plurigenera $P_m(X_t)$ is constant on $t \in S$ for any positive
integer $m$.}

\vskip 1pc

A morphism between complex varieties which is birationally equivalent to a 
projective morphism will be called an {\it algebraic morphism} in this paper.
The algebraicity assumption in the conjecture 
is slightly weaker than the projectivity.

This conjecture is confirmed by Iitaka [I1, I2] in the case in which 
$\text{dim }X_t = 2$ by using the classification theory of surfaces.
Nakayama [N1] proved that the conjecture follows if the minimal model exists 
for the family and the abundance conjecture holds for the generic fiber.
Thus the conjecture is true if $\text{dim }X_t = 3$ by [K4] and [KM].

On the other hand, Nakamura [Nm] provided a counterexample for the 
generalization of the conjecture in the case where the morphism $f$ is
not algebraic.  In his example, the central fiber $X_0$ is a quotient of a 
$3$-dimensional simply connected solvable Lie group by a discrete subgroup.
We note that $X_0$ is a non-K\"ahler manifold which 
has non-closed holomorphic $1$-forms.
So we only consider algebraic morphisms in this paper.
It is interesting to extend our results to the case in which 
the fibers are in Fujiki's class $\mathcal C$.

The following theorem of Siu was the starting point of the recent progress 
on this conjecture which we shall review.

\vskip 1pc

{\bf Theorem 1.2.} [Si]
{\em Let $S$ be a complex variety, and let
$f: \mathcal X \to S$ be a smooth projective morphism.
Assume that the generic fiber $X_{\eta}$ of $f$ is a variety of general type.
Then the plurigenera $P_m(X_t)$ is constant on $t \in S$ for any positive
integer $m$.}

\vskip 1pc

We have also a slightly stronger version:

\vskip 1pc

{\bf Theorem 1.2'.} [K5]
{\em Let $S$ be an algebraic variety, let $\mathcal X$ be a complex variety,
and let $f: \mathcal X \to S$ be a proper flat algebraic morphism.
Assume that the fibers $X_t = f^{-1}(t)$ have only
canonical singularities for any $t \in S$ and that the generic fiber
$X_{\eta}$ is a variety of general type.
Then the plurigenera $P_m(X_t)$ is constant on $t \in S$ for any positive
integer $m$.}

\vskip 1pc

According to Nakayama [N2], we define the {\it numerical Kodaira dimension} 
$\nu(X)$ as follows (this is $\kappa_{\sigma}(X)$ in [N2];
there is another version $\kappa_{\nu}(X)$ of numerical Kodaira dimension
in [N2] which we do not use).
Let $X$ be a compact complex manifold and let $k$ be a nonnegative integer.
We define $\nu(X) \ge k$ if there exist a divisor $H$ on $X$
and a positive number $c$
such that $\text{dim }H^0(X, mK_X + H) \ge cm^k$ for any sufficiently large 
and divisible $m$. If there is no such $k$, then we put $\nu(X) = - \infty$.
It is easy to see that $\kappa(X) \le \nu(X) \le \text{dim }X$.
By the Kodaira lemma, $\kappa(X) = \text{dim }X$ if and only if
$\nu(X) = \text{dim }X$.
The {\it abundance conjecture} states that the equality 
$\kappa(X) = \nu(X)$ always holds.
Nakayama confirmed this conjecture in the case $\nu(X) = 0$ ([N2]).

By considering $mK_X + H$ instead of $mK_X$, Nakayama obtained the following:

\vskip 1pc

{\bf Theorem 1.3.} [N3]
{\em Let $S$ be an algebraic variety, let $\mathcal X$ be a complex veriety,
and let $f: \mathcal X \to S$ be a proper flat algebraic morphism.
Assume that the fibers $X_t = f^{-1}(t)$ have only
canonical singularities for any $t \in S$.
Then the numerical Kodaira dimension $\nu(X_t)$ 
is constant on $t \in S$.
In particular, if one fiber $X_0$ is of general type, 
then so are all the fibers.}

\vskip 1pc

For the finer classification of algebraic varieties, 
it is useful to consider not only the discrete invariants $P_m(X)$ 
but also the infinite sum of vector spaces 
$$
R(X) = \bigoplus_{m \ge 0} H^0(X, mK_X)
$$ 
which has a natural graded ring structure over $\Bbb C = H^0(X, \mathcal O_X)$.
This continuous invariant $R(X)$, called the {\it canonical ring} of $X$,
is also independent of the birational models of $X$.
It is conjectured that R(X) is always finitely generated as a graded 
$\Bbb C$-algebra.  If this is the case, then 
$\text{Proj }R(X)$ is called a {\it canonical model} of $X$.

A {\it canonical singularity} (resp. {\it terminal singularity})
is defined as a singularity which may appear on 
a canonical model of a variety of general type whose canonical 
ring is finitely generated (resp. on a minimal model on an algebraic 
variety).  The formal definition by Reid is as follows:
a normal variety $X$ is said to have only 
canonical singularities (resp. terminal singularities) if the canonical 
divisor $K_X$ is $\Bbb Q$-Cartier and, for a resolution of singularities
$\mu: Y \to X$ which has exceptional divisors $F_j$, 
if we write $\mu^*K_X = K_Y + \sum_j a_j F_j$, then $a_j \le 0$
(resp. $a_j < 0$) for all $j$.

For example, the canonical singularities in dimension $2$
have been studied extensively.  They are called in many names such as
du Val singularities, rational double points, simple singularities, or
A-D-E singularities.
The terminal singularity in dimension $2$ is smooth,
and the terminal singularities in dimension $3$ are classified by
Mori and others (cf. [R2]).

Let us consider the subset of a Hilbert scheme
with a given Hilbert polynomial which consists of points corresponding to
the canonical models of varieties of general type.
This set should be open from the view point of 
the moduli problem of varieties (cf. [V2]).
The following is a local version of Theorem 1.2' and says that this is the case
(this result was previously known up to dimension $2$):

\vskip 1pc

{\bf Theorem 1.4.} [K5]
{\em Let $f: \mathcal X \to B$ be a flat morphism from a germ of an algebraic variety
to a germ of a smooth curve.
Assume that the central fiber $X_0 = f^{-1}(P)$ has only canonical 
singularities.
Then so has the total space $\mathcal X$ as well as any fiber $X_t$ of $f$.
Moreover, if $\mu: V \to \mathcal X$ is a birational morphism from a normal variety
with the strict transform $X$ of $X_0$, 
then $K_V + X \ge \mu^*(K_{\mathcal X} + X_0)$.}

\vskip 1pc

The folloing theorem answers a similar question for the deformations
of minimal models (this was previously known up to dimension $3$):

\vskip 1pc

{\bf Theorem 1.5.} [N3]
{\em Let $f: \mathcal X \to B$ be a flat morphism from a germ of an algebraic variety
to a germ of a smooth curve.
Assume that the central fiber $X_0 = f^{-1}(P)$ has only terminal
singularities.
Then so has the total space $\mathcal X$ as well as any fiber $X_t$ of $f$.
Moreover, if $\mu: V \to \mathcal X$ is a birational morphism from a normal variety
with the strict transform $X$ of $X_0$, 
then the support of $K_V + X - \mu^*(K_{\mathcal X} + X_0)$
contains all the exceptional divisors of $\mu$.}

\vskip 1pc

The following theorem, which is stronger than Theorem 1.2', says that 
only the abundance conjecture for the generic fiber implies the
deformation invariance of the plurigenera:

\vskip 1pc

{\bf Theorem 1.6.} [N3]
{\em Let $S$ be an algebraic variety, let $\mathcal X$ be a complex variety,
and let $f: \mathcal X \to S$ be a proper flat algebraic morphism.
Assume that the fibers $X_t = f^{-1}(t)$ have only
canonical singularities and that $\kappa(X_{\eta}) = \nu(X_{\eta})$ 
for the generic fiber $X_{\eta}$ of $f$.
Then the plurigenera $P_m(X_t)$ is constant on $t \in S$ for any positive
integer $m$.}

\vskip 1pc

Now we explain the idea of the proofs.
Since we assumed the algebraicity of varieties, there exist divisors
which are big. Hence we can use the vanishing theorems of Kodaira type as in 
[K1] and [V1] (cf Theorem 2.6).
Indeed, if $K_{X_0}$ is nef and big for the central fiber $X_0$
in Theorem 1.2, 
then the extendability of pluricanonical forms follows 
immediately from the vanishing theorem.

Thus the problem is to extract the nef part from the big divisor $K_{X_0}$.
This is similar to the {\it Zariski decomposition} problem (cf. [K3]):
Let $X$ be a smooth projective variety of general type.
If we fix a positive integer $m$, then there exists a projective 
birational morphism $\mu_m: Y_m \to X$ such that $\mu_m^*(mK_{X})$
is decomposed into the sum of the free part and the fixed part:
$\mu_m^*(mK_{X}) = P_m + M_m$.
If there exists one $\mu: Y \to X$ 
which serves as the $\mu_m$ simultaneously for all
$m$, then $P = \text{sup}_{m > 0}P_m/m$ is the desired nef 
part, and the decomposition $\mu^*K_X = P + N$ 
in $\text{Div}(Y) \otimes \Bbb R$
for $N = \text{inf}_{m > 0}N_m/m$ gives the
Zariski decomposition of $K_X$. 
The difficulty arises when we have an infinite tower of blow-ups.
It is known that if the Zariski decomposition of the canonical divisor 
exists, then the canonical ring $R(X)$ is finitely generated ([K3]).

So we use instead the concept of {\it multiplier ideal sheaf} which was first 
introduced by Nadel [Nd].
We consider the series of ideal sheaves on $X_0$ instead of the decompositions 
on the series of blow-ups.
Since the structure sheaf of $X_0$ is noetherian, we do not have the difficulty
of the infitity in this case; we take just the union of the ideals
(cf. Definitions 2.5 and 2.10).

The remaining thing to be proved is the compatibility of 
the multiplier ideal sheaves on $X_0$ and 
on the total space $\mathcal X$ constructed similarly for $K_{\mathcal X}$.
This is proved by a tricky induction on $m$ discovered by Siu (cf. Lemma 3.6).

The theorems in the introduction will be 
reduced to Theorems A, B and C in \S 2,
which will be proved by using vanishing theorems in \S 3.

We use the following terminology besides those in [KMM].
Let $f: X \to S$ be a morphism of algebraic varieties.
A sheaf $\mathcal F$ on $X$ is said to be {\it $f$-generated}
if the natural homomorphism $f^*f_*\mathcal F \to \mathcal F$
is surjective. 
A Cartier divisor $D$ on $X$ is called 
{\it $f$-effective} (resp. {\it $f$-free}) if 
$f_*\mathcal O_X(D) \ne 0$ 
(resp. $\mathcal O_X(D)$ is $f$-generated).
A $\Bbb Q$-Cartier divisor $D$ on $X$ is said to be 
{\it $f$-$\Bbb Q$-effective} (resp. {\it $f$-semi-ample}) if 
there exists a positive integer $m$ such that $mD$ is a $f$-effective 
(resp. $f$-free) Cartier divisor.
A $\Bbb Q$-Cartier divisor $D$ on $X$ is said to be 
{\it $f$-pseudo-effective} if $D + H$ is $f$-$\Bbb Q$-effective for any
$f$-ample $\Bbb Q$-Cartier divisor $H$. 

All varieties and morphisms are defined over the complex number field $\Bbb C$
in this paper.

%%%%%%%%%%%%%%%%%%%%%%%%%%%%%
%%%%%%%%%%%%%%%%%%%%%%%%%%%%%
%%%%%%%%%%%%%%%%%%%%%%%%%%%%%
\section{Main theorems}

{\bf Setup 2.1.}

{\em We fix the following notation in Theorems A, B and C below.

(1) $V$ is a smooth algebraic variety.

(2) $X$ is a smooth divisor on $V$.

(3) $S$ is a germ of an algebraic variety.

(4) $\pi: V \to S$ is a projective morphism:}
$$
\begin{CD}
X \subset V @>{\pi}>> S.
\end{CD}
$$
%\enddefinition

A divisor $D$ on $V$ will be called {\it $\pi$-effective for the pair} $(V, X)$
if the natural homomorphism
$\pi_*\mathcal O_V(D) \to \pi_*\mathcal O_X(D \vert_X)$ is not zero.
It is called {\it $\pi$-$\Bbb Q$-effective for the pair} $(V, X)$
if $mD$ is $\pi$-effective for the pair $(V, X)$ for some positive integer $m$.
$D$ is said to be {\it $\pi$-big for the pair} $(V, X)$ if 
we can write $mD = A + B$ for a positive integer $m$, $\pi$-ample 
divisor $A$ and a $\pi$-effective divisor $B$ for the pair $(V, X)$.

\vskip 1pc

{\bf 2.2. Theorem A.} 
{\em Assume that $K_V + X$ is $\pi$-big for the pair $(V, X)$.
Then the natural homomorphism
$\pi_*\mathcal O_V(m(K_V + X)) \to \pi_*\mathcal O_X(mK_X)$ is surjective 
for any positive integer $m$.}

\vskip 1pc

\begin{proof}[Theorem A implies Theorem 1.4]
By the resolution of singularities, we construct a projective 
birational morphism $\mu: V \to \mathcal X$ from a smooth variety
such that the strict transform $X$ of the central fiber $X_0$ is also
smooth. We take $S = \mathcal X$ and $\pi = \mu$.
Since $\pi$ and $\pi \vert_X$ are birational morphisms, we can see that
$K_V + X$ is $\pi$-big for the pair $(V, X)$.

Let $m$ be a positive integer such that $mK_{X_0}$ is a Cartier divisor.
Since $X_0$ has only canonical singularities, a nowhere vanishing section $s_0$
of $\mathcal O_{X_0}(mK_{X_0})$ lifts to a section of $\mathcal O_X(mK_X)$, 
which in turn extends to a section of $\mathcal O_V(m(K_V + X))$ by Theorem A.  
Therefore, $s_0$ extends to a nowhere vanishing 
section of $\mathcal O_{\mathcal X}(m(K_{\mathcal X} + X_0))$ which lifts to a section of
$\mathcal O_V(m(K_V + X))$.
Hence $K_{\mathcal X} + X_0$ is a $\Bbb Q$-Cartier divisor and 
$K_V + X \ge \mu^*(K_{\mathcal X} + X_0)$ in this case.
Since any $\mu$ is dominated by some $\mu$ as above, we have the
last assertion.
Since $\mu^*X_0 \ge X$, it follows that 
$\mathcal X$ has only canonical singularities.
\end{proof}

\begin{proof}[Theorem A implies Theorem 1.2']
We may assume that $S$ is a germ of a smooth curve.
By the resolution of singularities, we construct a proper birational morphism
$\mu: V \to \mathcal X$ from a smooth variety
such that the strict transform $X$ of the central fiber $X_0$ is also
smooth and that $\pi = f \circ \mu$ is a projective morphism.

Let $A$ be a $\pi$-very ample divisor on $V$.
Since $K_{\mathcal X}$ is $f$-big, there exists a positive integer $m_1$ such 
that $f_*\mathcal O_{\mathcal X}(m_1K_{\mathcal X} - \mu_*A) \ne 0$, although $\mu_*A$
is a Weil divisor which may not be a Cartier divisor.
Therefore, there exists an effective Weil divisor $\bar B$ on $\mathcal X$ whose 
support does not contain $X_0$ such that $m_1K_{\mathcal X} \sim \mu_*A + \bar B$.
Then $K_V + X$ is $\pi$-big for the pair $(V, X)$ by the last assertion of 
Theorem 1.4.
Since $\mathcal X$ has only canonical singularities by Theorem 1.4,  
Theorem A implies that the natural homomorphism
$f_*\mathcal O_{\mathcal X}(mK_{\mathcal X}) \to H^0(X_0, mK_{X_0})$ is 
surjective for any positive integer $m$.
\end{proof}

\vskip 1pc

{\bf 2.3. Theorem B.} 
{\em Let $H$ be a $\pi$-very ample divisor on $V$.
Assume that $K_X$ is $\pi$-pseudo-effective.
Then $K_V + X$ is also $\pi$-pseudo-effective, and the natural homomorphism
$\pi_*\mathcal O_V(m(K_V + X) + H) \to \pi_*\mathcal O_X(mK_X + H \vert_X)$ 
is surjective for any positive integer $m$.}

\vskip 1pc

\begin{proof}[Theorem B implies Theorem 1.3]
We may assume that $S$ is a germ of a smooth curve.
We define $V, X$ and $\pi$ as in the proof that Theorem A implies Theorem 1.2'.
If $\nu(X_0) = - \infty$, then $K_{X_0}$ is not pseudo-effective, 
and $\nu(X_{\eta}) = - \infty$ for the generic fiber $X_{\eta}$
by the upper semi-continuity theorem.
Otherwise, the rest of the proof is 
similar to the proof that Theorem A implies Theorem 1.2'.
\end{proof}

\begin{proof}[Theorem B implies Theorem 1.5]
There exist a resolution of singularities $\mu: V \to \mathcal X$ and 
an effective divisor $E$ which is supported on the exceptional 
locus of $\mu$ such that $H = - E$ is $\mu$-very ample.
Let $S = \mathcal X$ and $\pi = \mu$.
Let $X$ be the strict transform of $X_0$ which 
is assumed to be smooth.
Since $X_0$ has only terminal singularities, there exists a positive 
integer $m$ such that $mK_{X_0}$ is a Cartier divisor and that
$mK_X - m\mu^*K_{X_0} \ge E \vert_X$.
Thus a nowhere vanishing section $s_0$
of $\mathcal O_{X_0}(mK_{X_0})$ lifts to a section of $\mathcal O_X(mK_X - E \vert_X)$,
which in turn extends to a section of $\mathcal O_V(m(K_V + X) - E)$ by Theorem B.
Therefore, $s_0$ extends to a nowhere vanishing 
section of $\mathcal O_{\mathcal X}(m(K_{\mathcal X} + X_0))$ which lifts to a section of
$\mathcal O_V(m(K_V + X) - E)$.
Hence $K_V + X - \mu^*(K_{\mathcal X} + X_0) \ge \frac 1m E$ in this case.
Since any $\mu$ is dominated by some $\mu$ as above, we have the
last assertion.
Since $\mu^*X_0 \ge X$, it follows that 
$\mathcal X$ has only terminal singularities.
\end{proof}

\vskip 1pc

{\bf 2.4. Theorem C.} 
{\em Let $X_{\xi}$ and $V_{\eta}$ be the generic fibers of $\pi: X \to \pi(X)$ and
$\pi: V \to \pi(V)$, respectively.
Assume that $K_V + X$ is $\pi$-$\Bbb Q$-effective for the pair $(V, X)$,
$\text{dim }X_{\xi} = \text{dim }V_{\eta}$, 
and that $\nu(X_{\xi}) = \nu(V_{\eta}) = \kappa(V_{\eta})$. 
Then the natural homomorphism
$\pi_*\mathcal O_V(m(K_V + X)) \to \pi_*\mathcal O_X(mK_X)$ is surjective 
for any positive integer $m$.}

\vskip 1pc

\begin{proof}[Theorem C implies Theorem 1.6]
We may assume that $S$ is a germ of a smooth curve.
We define $V, X$ and $\pi$ as in the proof that Theorem A implies Theorem 1.2'.
We may assume that $K_{X_0}$ is pseudo-effective.
By Theorem 1.3, we have $\nu(X_{\eta}) = \nu(X_0) \ge 0$, where
$X_{\eta}$ is the generic fiber of $f$.
By the assumption on the abundance, it follows that $K_V + X$ becomes 
$\pi$-$\Bbb Q$-effective for the pair $(V, X)$.
The rest is similar to the proof that Theorem A implies Theorem 1.2'.
\end{proof}

\vskip 1pc

{\bf Definition 2.5.}
Let $X$ be a smooth complex variety and $D$ a divisor on $X$.
Let $\mu: Y \to X$ be a proper birational morphism from a smooth variety $Y$.
Assume that there exists a decomposition $\mu^*D = P + M$ 
in $\text{Div}(Y) \otimes \Bbb R$ such that $P$ is $\mu$-nef and 
$M$ is effective having a normal crossing support.
The {\it multiplier ideal sheaf} $\mathcal I_M$ is defined by the following 
formula:
$$
\mu_*\mathcal O_Y(\ulcorner P \urcorner + K_Y) 
= \mathcal I_M(D + K_X).
$$
We note that $\mathcal I_M$ is a coherent sheaf of ideals of $\mathcal O_X$ which is 
determined only by $M$ and $\mu$.
If $\nu: Y' \to Y$ is another proper birational morphism, then it is easy to
see that $\mathcal I_M = \mathcal I_{\nu^*M}$.

\vskip 1pc

The following vanishing theorem of Kawamata-Viehweg type ([KMM, 1.2.3]) 
is the main tool for the proof of Theorems A and B.  
Nadel's vanishing theorem [Nd] and
Ohsawa-Takegoshi's extension theorem [OT] played the same role in Siu's 
proof.

\vskip 1pc

{\bf Theorem 2.6.} 
{\em Let $f: X \to S$ be a proper algebraic morphism from a smooth complex manifold 
to an algebraic variety. Let $L$ be an $\Bbb R$-divisor on $X$ which is 
$f$-nef and $f$-big, and whose fractional part has a normal crossing 
support. Then 
$$
R^pf_*\mathcal O_X(\ulcorner L \urcorner + K_X) = 0
$$
for any positive integer $p$.}

\vskip 1pc

{\bf Corollary 2.7.}
{\em In the situation of Definition 2.5, 
let $S$ be an algebraic variety, and let 
$f: X \to S$ be a proper algebraic morphism.
If $P$ is $f \circ \mu$-nef and $f \circ \mu$-big, then
$$
R^pf_*\mathcal I_M(D + K_X) = 0
$$
for any positive integer $p$.}

\vskip 1pc

\begin{proof}
We have $R^p(f \circ \mu)_*\mathcal O_Y(\ulcorner P \urcorner + K_Y) = 0$ and
$R^p\mu_*\mathcal O_Y(\ulcorner P \urcorner + K_Y) = 0$ for any positive integer
$p$.
\end{proof}

We need a vanishing theorem of Koll\'ar type ([Ko]) in order to prove 
Theorem C:

\vskip 1pc

{\bf Theorem 2.8.}
{\em Let $f: X \to S$ be a proper algebraic morphism from a smooth complex manifold 
to an algebraic variety. Let $L$ be a $\Bbb Q$-divisor on $X$ which is 
$f$-semi-ample and whose fractional part has a normal crossing support.
Let $D$ be an effective divisor on $X$ such that
$mL - D$ is $f$-effective for a positive integer $m$.
Then a natural homomorphism
$$
R^pf_*\mathcal O_X(\ulcorner L \urcorner + K_X)
\to R^pf_*\mathcal O_X(\ulcorner L \urcorner + D + K_X)
$$
is injective for any nonnegative integer $p$.}

\vskip 1pc

\begin{proof}
We may assume that $S$ is affine and $f$ is projective by Theorem 2.6.
By compactifying $S$, adding the pull-back of an ample divisor of $S$ to 
$L$, and using the Serre vanishing theorem, 
we reduce the assertion to the injectivity of the homomorphism
$$
H^p(X, \mathcal O_X(\ulcorner L \urcorner + K_X))
\to H^p(X, \mathcal O_X(\ulcorner L \urcorner + D + K_X))
$$
in the case in which $X$ is projective and $L$ is semi-ample.
This is just [K2, Theorem 3.2].
\end{proof}

\vskip 1pc

{\bf Corollary 2.9.}
{\em In the situation of Corollary 2.7,
if $f$ is surjective and $P$ is $f \circ \mu$-semi-ample, then
$$
R^pf_*\mathcal I_M(D + K_X)
$$
is torsion free for any nonnegative integer $p$.}

\vskip 1pc

{\bf Definition 2.10.}
In the situation of Setup 2.1, 
we define several kinds of multiplier ideal sheaves 
$\mathcal J^0_D, \mathcal I^0_D, \mathcal J^1_D, \mathcal I^1_D$ on $X$ in the following.

Let $D$ be a divisor on $X$ which is $\pi$-$\Bbb Q$-effective.
For each positive integer $m$ such that $mD$ is $\pi$-effective, 
we construct a proper birational morphism $\mu_m: Y_m \to X$
from a smooth variety such that the following conditions are satisfied:
there is a decomposition $\mu_m^*(mD) = P_m + M_m$ in $\text{Div}(Y_m)$
such that $P_m$ is $\pi \circ \mu_m$-free,
$M_m$ is effective and has a normal crossing support, and that
the natural homomorphism
$(\pi \circ \mu_m)_*\mathcal O_{Y_m}(P_m) \to \pi_*\mathcal O_X(mD)$ 
is an isomorphism.
We define
$$
\mathcal J^0_D = \bigcup_m \mathcal I_{\frac 1m M_m}
$$
where the union is taken for all positive integers $m$ such that $mD$ is 
$\pi$-effective.
Since $X$ is noetherian, there exists a positive integer $m$ such that
$\mathcal J^0_D = \mathcal I_{\frac 1m N_m}$.

In the case in which $D$ itself is $\pi$-effective,
we define $\mathcal I^0_D = \mathcal I_{M_1}$. 
We have $\mathcal I^0_D \subset \mathcal J^0_D$.

Let $D$ be a divisor on $V$ which is $\pi$-$\Bbb Q$-effective 
for the pair $(V, X)$.
For each positive integer $m$ such that $mD$ is $\pi$-effective for the
pair $(V, X)$, 
we construct a proper birational morphism $\mu_m: W_m \to V$
from a smooth variety with the strict transform $Y_m$ of $X$ in
$W_m$ being smooth and such that the following conditions are satisfied:
there is a decomposition $\mu_m^*(mD) = Q_m + N_m$
in $\text{Div}(W_m)$ such that $Q_m$ is $\pi \circ \mu_m$-free,
$N_m$ is effective, $N_m + Y_m$ has a normal crossing support, and that
the natural homomorphism
$(\pi \circ \mu_m)_*\mathcal O_{W_m}(Q_m) \to \pi_*\mathcal O_V(mD)$ is an isomorphism.
Since $mD$ is $\pi$-effective for the
pair $(V, X)$, $Y_m$ is not contained in the support of $N_m$.
We define 
$$
\mathcal J^1_D = \bigcup_m \mathcal I_{\frac 1m N_m \vert_{Y_m}}
$$
where the union is taken for all positive integers $m$ such that $mD$ is 
$\pi$-effective for the pair $(V, X)$.
Since $X$ is noetherian, there exists a positive integer $m$ such that
$\mathcal J^1_D = \mathcal I_{\frac 1m N_m \vert_{Y_m}}$.

In the case in which $D$ itself is $\pi$-effective for the pair $(V, X)$,
we define $\mathcal I^1_D = \mathcal I_{N_1 \vert_{Y_1}}$. 
We have $\mathcal I^1_D \subset \mathcal J^1_D$.

If we define $M_m$ for $D \vert_X$ as before, then we have
$N_m \vert_{Y_m} \ge M_m$ as divisors on $Y_m$.
Hence $\mathcal J^1_D \subset \mathcal J^0_{D \vert_X}$.
We also have $\mathcal I^1_D \subset \mathcal I^0_{D \vert_X}$ if $D$ is 
$\pi$-effective for the pair $(V, X)$.

\vskip 1pc

In the proof of the main results in the next section, 
the question on the extendability of global sections of the sheaves
$\mathcal O_X(D \vert_X)$ to those of $\mathcal O_V(D)$ 
will be reduced to the comparison
of the multiplier ideal sheaves such as $\mathcal J^0_{D \vert_X}$ 
and $\mathcal J^1_D$.

The following lemma enables us to deduce
the inclusion of the sheaves from
the inclusion of their direct image sheaves.
A theorem of Skoda [Sk] was used in [Si] to prove the corresponding statement
in the analytic setting.

\vskip 1pc

{\bf Lemma 2.11.} ([N3])
{\em Let $f: X \to S$ be a projective morphism of algebraic schemes such that 
$n = \text{dim }X$.
Let $\mathcal H$ be an $f$-very ample invertible sheaves on $X$, and 
let $\mathcal F$ be a coherent sheaf on $X$.
Assume that $R^pf_*(\mathcal F \otimes \mathcal H^{\otimes m}) = 0$ 
for any $p > 0$ and $m \ge 0$.
Then the sheaf $\mathcal F \otimes \mathcal H^{\otimes n}$ is $f$-generated.}

\vskip 1pc

\begin{proof}
We proceed by induction on $n$.
Let us take an arbitrary point $x \in X$.
We may assume that $H^0_{\{x\}}(\mathcal F) = 0$ if we replace $\mathcal F$ by
its quotient by the torsion subsheaf supported at $x$.
Let $X'$ be a general member in the linear system $\vert \mathcal H \vert$ 
passing through $x$.
Let $\mathcal H' = \mathcal H \otimes \mathcal O_{X'}$
and $\mathcal F' = \mathcal F \otimes \mathcal H'$.
We have an exact sequence
$$
0 \to \mathcal F \to \mathcal F \otimes \mathcal H \to \mathcal F' \to 0.
$$
By the induction hypothesis, $\mathcal F' \otimes \mathcal H^{\prime \otimes (n-1)}$ 
is $f$-generated. We have $R^1f_*(\mathcal F \otimes \mathcal H^{\otimes (n-1)}) = 0$, 
hence our assertion is proved.
\end{proof}

\vskip 1pc

{\bf Remark 2.12.}
We note that Theorems A, B and C also hold in the case in which 
$X$ is reducible.
So we can apply them when $f: \mathcal X \to S$ has reducible fibers 
in Theorems 1.2', 1.3 and 1.6 ([N3]).

%%%%%%%%%%%%%%%%%%%%%%%%%%%%%
%%%%%%%%%%%%%%%%%%%%%%%%%%%%%
%%%%%%%%%%%%%%%%%%%%%%%%%%%%%
\section{Proof of the main theorems}

{\bf Notation 3.1.}
{\em We fix divisors $H$ and $A$ on $V$ such that 
$H$ is $\pi$-very ample and 
$A = (\text{dim }X + 1)H$.}

\vskip 1pc

{\bf Lemma 3.2.}
{\em Let $D_1$ and $D$ be $\pi$-effective and $\pi$-$\Bbb Q$-effective divisors 
on $X$ (resp. $\pi$-effective and $\pi$-$\Bbb Q$-effective divisors on $V$ 
for the pair $(V, X)$), respectively. 
Then the sheaves $\mathcal I^0_{D_1}(D_1 + A \vert_X + K_X)$ and
$\mathcal J^0_D(D + A \vert_X + K_X)$ 
(resp. $\mathcal I^1_{D_1}(D_1 \vert_X + A \vert_X + K_X)$ and
$\mathcal J^1_D(D \vert_X + A \vert_X + K_X)$) are $\pi$-generated.}

\vskip 1pc

\begin{proof}
Since $A - (\text{dim }X)H$ is $\pi$-ample,
we apply Corollary 2.7 and Lemma 2.11.
\end{proof}

\vskip 1pc

{\bf Lemma 3.3.} 
{\em Let $D_1$ and $D$ be $\pi$-effective and $\pi$-$\Bbb Q$-effective divisors 
on $X$ (resp. $\pi$-effective and $\pi$-$\Bbb Q$-effective divisors 
on $V$ for the pair $(V, X)$). Then

(1) $\mathcal J^i_{\alpha D} \subset \mathcal J^i_D$ for $i = 0,1$
if $\alpha \in \Bbb Q$ and $\alpha > 1$.

(2) $\mathcal I^i_{D_1} \subset \mathcal I^i_{D_1 + L}$ and 
$\mathcal J^i_D \subset \mathcal J^i_{D + L}$ for $i = 0,1$
if $L$ is a $\pi$-free divisor on $X$ (resp. $V$).

(3) $\text{Im}(\pi_*\mathcal O_V(D_1) \to \pi_*\mathcal O_X(D_1 \vert_X))
\subset \pi_*\mathcal I^1_{D_1}(D_1 \vert_X)$ (in the latter case only).}

\vskip 1pc

\begin{proof}
(1) and (2) are clear.
(3) follows from $K_{Y_1} \ge \mu_1^*K_X$.
\end{proof}

The vanishing theorem is used to prove the following lemma.

\vskip 1pc

{\bf Lemma 3.4.} 
{\em Let $D$ be a $\pi$-$\Bbb Q$-effective divisor
for the pair $(V, X)$. Then}
$$
\pi_*\mathcal J^1_D(D \vert_X + K_X) 
\subset \text{Im}(\pi_*\mathcal O_V(D + K_V + X) \to 
\pi_*\mathcal O_X(D \vert_X + K_X)).
$$

\vskip 1pc

\begin{proof}
We have an exact sequence
\[
\begin{split}
&0 \to \mathcal O_{W_m}(\ulcorner \frac 1m Q_m \urcorner + K_{W_m}) 
\to \mathcal O_{W_m}(\ulcorner \frac 1m Q_m \urcorner + K_{W_m} + Y_m) \\
&\to \mathcal O_{Y_m}(\ulcorner \frac 1m Q_m \vert_{Y_m} \urcorner + K_{Y_m}) 
\to 0.
\end{split}
\]
If $D$ is $\pi$-big for the pair $(V, X)$, then we have 
$R^1(\pi \circ \mu_m)_*\mathcal O_{W_m}(\ulcorner \frac 1m Q_m \urcorner 
+ K_{W_m}) = 0$ by Theorem 2.6.
In the general case, since $Q_m$ is $\pi \circ \mu_m$-free, the homomorphism
$$
R^1(\pi \circ \mu_m)_*\mathcal O_{W_m}(\ulcorner \frac 1m Q_m \urcorner 
+ K_{W_m}) \to
R^1(\pi \circ \mu_m)_*\mathcal O_{W_m}(\ulcorner \frac 1m Q_m \urcorner 
+ K_{W_m} + Y_m)
$$
is injective by Theorem 2.8.
Anyway, we have a surjective homomorphism
\[
\begin{split}
&\pi_*\mathcal O_V(D + K_V + X) \supset 
(\pi \circ \mu_m)_*\mathcal O_{W_m}(\ulcorner \frac 1m Q_m \urcorner 
+ K_{W_m} + Y_m) \\
&\twoheadrightarrow (\pi \circ \mu_m)_*\mathcal O_{Y_m}
(\ulcorner \frac 1m Q_m \vert_{Y_m} \urcorner + K_{Y_m})
= \pi_*\mathcal J^1_D(D \vert_X + K_X),
\end{split}
\]
hence the assertion. 
\end{proof}

\vskip 1pc

{\bf Corollary 3.5.}
{\em (1) Let $D$ be a $\pi$-$\Bbb Q$-effective divisor
for the pair $(V, X)$. Then $D + A + K_V + X$ is $\pi$-effective
for the pair $(V, X)$.

(2) If $D_1$ and $D_1 + K_V + X$ are $\pi$-effective for the pair $(V, X)$, 
and $D$ and $D + K_V + X$ are $\pi$-$\Bbb Q$-effective
for the pair $(V, X)$.  Then}
\[
\begin{split}
\pi_*\mathcal J^1_{D_1}(D_1 \vert_X + K_X) 
&\subset \pi_*\mathcal I^1_{D_1 + K_V + X}(D_1 \vert_X + K_X) \\
\pi_*\mathcal J^1_D(D \vert_X + K_X) 
&\subset \pi_*\mathcal J^1_{D + K_V + X}(D \vert_X + K_X).
\end{split}
\]

\vskip 1pc

\begin{proof}
(1) By Lemma 3.2, we have
$\pi_*\mathcal J^1_D(D \vert_X + A \vert_X + K_X) \ne 0$.
By Lemma 3.3 (2), $\pi_*\mathcal J^1_{D + A}(D \vert_X + A \vert_X + K_X) \ne 0$.
Then by Lemma 3.4, we obtain our assertion.
(2) follows from Lemmas 3.3 (3) and 3.4.
\end{proof}

The following is the key lemma for Theorems A and C.

\vskip 1pc

{\bf Lemma 3.6.} 
{\em Assume that $K_V + X$ is $\pi$-$\Bbb Q$-effective for the pair $(V, X)$. Then
$$
\mathcal J^0_{mK_X} \subset \mathcal I^1_{m(K_V + X) + A}
$$
for any non-negative integer $m$. }

\vskip 1pc

\begin{proof}
We note that $m(K_V + X) + A$ is $\pi$-effective for any non-negative
integer $m$ by Corollary 3.5 (1).
We proceed by induction on $m$.
If $m = 0$, then the assertion is obvious, because $A$ is $\pi$-free.
Assume that the assertion is true for $m - 1$.
By the induction hypothesis and Corollary 3.5 (2),
\[
\begin{split}
&\pi_*\mathcal J^0_{(m-1)K_X}(mK_X + A \vert_X) 
\subset \pi_*\mathcal I^1_{(m-1)(K_V + X)+ A}(mK_X + A \vert_X) \\
&\subset \pi_*\mathcal I^1_{m(K_V + X) + A}(mK_X + A \vert_X).
\end{split}
\]
Since $\pi_*\mathcal J^0_{(m-1)K_X}(mK_X + A \vert_X)$ is $\pi$-generated
by Lemma 3.2, it follows that 
$$
\mathcal J^0_{(m-1)K_X} \subset \mathcal I^1_{m(K_V + X) + A}.
$$
Since $\mathcal J^0_{mK_X} \subset \mathcal J^0_{(m-1)K_X}$
by Lemma 3.3 (1), we are done.
\end{proof}

\begin{proof}[Proof of Theorem A]
Since $K_V + X$ is $\pi$-big for the pair $(V, X)$, 
there exists a positive integer $m_0$ such that $m_0(K_V+X) \sim A + B$ for 
an effective divisor $B$ whose support does not contain $X$.
By Lemmas 3.3 (1) and 3.6, we have
\[
\begin{split}
&\mathcal J^0_{mK_X}(- B \vert_X) 
\subset \mathcal J^1_{m(K_V + X) + A}(- B \vert_X) \\
&\subset \mathcal J^1_{m(K_V + X) + A + B}
= \mathcal J^1_{(m + m_0)(K_V + X)} \subset \mathcal J^1_{m(K_V + X)}.
\end{split}
\]
This implies the following: for any positive integer $m$, 
there exists a positive integer $d$ such that,
if $Y$ is any smooth model of $X$ which dominates $Y_m$ and $Y_{dm}$ by
$\mu: Y \to X$, $\nu_m: Y \to Y_m$ and $\nu_{dm}: Y \to Y_{dm}$, then 
$$
- \nu_m^*M_m - \mu^*B \vert_X \le 
\ulcorner - \frac 1d \nu_{dm}^*N_{dm}\vert_{Y_{dm}} \urcorner 
+ K_Y - \mu^*K_X.
$$
If $m = en$ for positive integers $e$ and $n$ and if $Y$ dominates $Y_n$ by
a morphism $\nu_n: Y \to Y_n$, then we have
$$
- \nu_n^*M_n \le - \frac 1e \nu_m^*M_m
\le \frac 1e \mu^*B \vert_X + \frac 1e C_{m,d}
- \frac 1{de} \nu_{dm}^*N_{dm} \vert_{Y_{dm}} + \frac 1e (K_Y - \mu^*K_X)
$$
where $C_{m,d} = \ulcorner - \frac 1d \nu_{dm}^*N_{dm} \vert_{Y_{dm}} 
\urcorner + \frac 1d \nu_{dm}^*N_{dm} \vert_{Y_{dm}}$.
We take a large enough integer $e$ such that $(X, \frac 2e B \vert_X)$ 
is log terminal.
Then we have
$\llcorner \frac 2e \mu^*B \vert_X - (K_Y - \mu^*K_X) \lrcorner \le 0$.
Hence
$$
- \nu_n^*M_n 
\le \ulcorner - \frac 1{de}\nu_{dm}^*N_{dm} \vert_{Y_{dm}} \urcorner
+ K_Y - \mu^*K_X.
$$
Thus $\pi_*\mathcal O_X(nK_X) \subset \pi_*\mathcal J^1_{n(K_V + X)}(nK_X)
\subset \pi_*\mathcal J^1_{(n-1)(K_V + X)}(nK_X)$.
Therefore, we obtain our assertion by Lemma 3.4.
\end{proof}

We modify Lemma 3.6 for Theorem B as follows:

\vskip 1pc

{\bf Lemma 3.7.} 
{\em Assume that $K_X$ is $\pi$-pseudo-effective. 
Then $K_V + X$ is $\pi$-pseudo-effective,
and
$$
\mathcal J^0_{mK_X+ eH \vert_X} \subset \mathcal I^1_{m(K_V + X) + eH + A}.
$$ 
for any non-negative integer $m$ and any positive integer $e$.}

\vskip 1pc

\begin{proof}
Since $K_X$ is $\pi$-pseudo-effective, $mK_X + eH \vert_X$ is $\pi$-big 
for any non-negative integer $m$ and any positive integer $e$.
Thus the left hand side of the formula is well defined.
We shall prove that $m(K_V + X) + H + A$ is $\pi$-effective for the
pair $(V, X)$ for any non-negative integer $m$ in order to prove that
$K_V + X$ is $\pi$-pseudo-effective,
as well as the inclusion
$\mathcal J^0_{mK_X + H \vert_X} \subset \mathcal I^1_{m(K_V + X) + H + A}$
in the case $e = 1$ by induction on $m$.
The inclusion for general $e$ is proved similarly for each fixed $e$
by induction on $m$.

If $m = 0$, then the assertion is obvious.
Assume that the assertion is true for $m - 1$.
By the induction hypothesis and Lemma 3.4, we have
\[
\begin{split}
&\pi_*\mathcal J^0_{(m-1)K_X+ H \vert_X}(mK_X + H \vert_X + A \vert_X) \\
&\subset \pi_*\mathcal I^1_{(m-1)(K_V + X) + H + A}
(mK_X + H \vert_X + A \vert_X) \\
&\subset \text{Im}(\pi_*\mathcal O_V(m(K_V + X) + H + A) \to 
\pi_*\mathcal O_X(mK_X + H \vert_X + A \vert_X)).
\end{split}
\]
Since $\mathcal J^0_{(m-1)K_X + H \vert_X}(mK_X + H \vert_X + A \vert_X)$
is $\pi$-generated by Lemma 3.2, it follows that
$m(K_V + X) + H + A$ is $\pi$-effective for the pair $(V, X)$.
Then by Lemma 3.3 (3),
\[
\begin{split}
&\text{Im}(\pi_*\mathcal O_V(m(K_V + X) + H + A) \to 
\pi_*\mathcal O_X(mK_X + H \vert_X + A \vert_X)) \\
&\subset \pi_*\mathcal I^1_{m(K_V + X) + H + A}
(mK_X + H \vert_X + A \vert_X).
\end{split}
\]
Hence 
$$
\mathcal J^0_{(m-1)K_X + H \vert_X} 
\subset \mathcal I^1_{m(K_V + X) + H + A}.
$$
Since $\mathcal J^0_{mK_X + H \vert_X} \subset \mathcal J^0_{(m-1)K_X + H \vert_X}$
by Lemma 3.3 (1) and (2), we are done.
\end{proof}

\begin{proof}[Proof of Theorem B]
We fix $m$.  Since $m(K_V + X) + H$ is $\pi$-big for the pair $(V, X)$
by Lemma 3.7,
there exists a positive integer $m_0$ such that
$m_0(m(K_V + X) + H) \sim H + A + B$ for an effective divisor $B$
whose support does not contain $X$.
By Lemmas 3.3 (1), (2) and 3.7, we have
\[
\begin{split}
&\mathcal J^0_{mK_X+ eH \vert_X}(- B \vert_X) 
\subset \mathcal J^1_{m(K_V + X) + eH + A}(- B \vert_X) \\
&\subset \mathcal J^1_{m(K_V + X) + eH + A + B}
\subset \mathcal J^1_{(m_0 + 1)(m(K_V + X) + eH)} \\
&\subset \mathcal J^1_{m(K_V + X) + eH}.
\end{split}
\]
We have proper birational morphisms $\mu_{m,e}: W_{m,e} \to V$
from a smooth variety, a smooth strict transform $Y_{m,e}$ of $X$, 
and decompositions $\mu_{m,e}^*(mK_X + eH \vert_X)
= P_{m,e} + M_{m,e}$ and $\mu_{m,e}^*(m(K_V+X)+eH)
= Q_{m,e} + N_{m,e}$ as in Definition 2.10.

It follows that for any positive integers $m$ and $e$, 
there exists a positive integer $d$ such that,
if $Y$ is any smooth model of $X$ which dominates $Y_{m,e}$ and 
$Y_{dm,de}$ by $\mu: Y \to X$, $\nu_{m,e}: Y \to Y_{m,e}$ 
and $\nu_{dm,de}: Y \to Y_{dm,de}$, then 
$$
- \nu_{m,e}^*M_{m,e} - \mu^*B \vert_X \le 
\ulcorner - \frac 1d \nu_{dm,de}^*N_{dm,de}\vert_{Y_{dm,de}} 
\urcorner + K_Y - \mu^*K_X.
$$
If $m = en$ for a positive integer $n$, then we have
\[
\begin{split}
&- \nu_{n,1}^*M_{n,1} \le - \frac 1e \nu_{m,e}^*M_{m,e} \\
&\le \frac 1e \mu^*B \vert_X + \frac 1e C_{m,d,e}
- \frac 1{de} \nu_{dm,de}^*N_{dm,de} \vert_{Y_{dm,de}} 
+ \frac 1e (K_Y - \mu^*K_X)
\end{split}
\]
where $C_{m,d,e} = \ulcorner - \frac 1d \nu_{dm,de}^*N_{dm,de} 
\vert_{Y_{dm,de}} \urcorner + \frac 1d \nu_{dm,de}^*N_{dm,de} 
\vert_{Y_{dm,de}}$.
We take a large enough integer $e$ such that $(X, \frac 2e B \vert_X)$ 
is log terminal. Then 
$\llcorner \frac 2e \mu^*B \vert_X - (K_Y - \mu^*K_X) \lrcorner \le 0$,
hence
$$
- \nu_{n,1}^*M_{n,1} 
\le \ulcorner - \frac 1{de}\nu_{dm,de}^*N_{dm,de} \vert_{Y_{dm,de}} \urcorner
+ K_Y - \mu^*K_X.
$$
Therefore,
$$\pi_*\mathcal O_X(nK_X + H \vert_X) \subset 
\pi_*\mathcal J^1_{n(K_V + X) + H}(nK_X + H \vert_X),
$$
and the whole of $\pi_*\mathcal O_X(nK_X + H \vert_X)$ is extendable to $V$.
\end{proof}

\begin{proof}[Proof of Theorem C]
By the flattening and the normalization,
we construct a proper birational morphism $\mu: V' \to V$ 
from a normal variety $V'$, a smooth variety $\bar V$ with a structure
morphism $\beta: \bar V \to S$, and 
an equidimensional projective morphism 
$\alpha: V' \to \bar V$ which is birationally equivalent to the Iitaka 
fibration of $V$ over $S$.  We set $\pi' = \beta \circ \alpha
= \pi \circ \mu$:
$$
\begin{CD}
V @<{\mu}<< V' \\
@V{\pi}VV   @VV{\alpha}V \\
S @<{\beta}<< \bar V.
\end{CD}
$$

There exists a positive integer $m_1$ such that 
$m_1\mu^*(K_V + X) = \alpha^*\bar L + E$ 
for a $\beta$-big divisor $\bar L$ on $\bar V$ and
an effective Cartier divisor $E$ on $V'$.
Since the Kodaira dimension of the general fiber of $\alpha$ is zero and 
$\alpha$ is equidimensional, we may assume that the natural homomorphisms
$\beta_*\mathcal O_{\bar V}(m\bar L) \to \pi_*\mathcal O_V(mm_1(K_V + X))$ 
are bijective for any positive integer $m$.
This means in particular that for any prime divisor of $\bar V$, 
there exists a prime divisor of $V'$ lying above it which is not contained in 
the support of $E$.

Since $\kappa(V_{\eta}) = \nu(V_{\eta})$, 
the numerical Kodaira dimension of the generic
fiber of $\alpha$ is zero by [N2, 7.4.3].
Hence there exists a positive integer $m_2$ such that 
the sheaf $\mathcal F = \alpha_*\mathcal O_{V'}(\mu^*A + mE)$ on $\bar V$
is independent of the integer $m$ if $m \ge m_2$.

Let $X'$ be the strict transform of $X$ by $\mu$.
Since $\text{dim }X_{\xi} = \text{dim }V_{\eta}$, 
we have $\alpha(X') = \bar X \ne \bar V$.
We may assume that $\bar X$ is a smooth divisor on $\bar V$.
Since $\mathcal F$ is torsion free,
there exists a $\beta$-very ample divisor $\bar A$ on $\bar V$ 
such that there are injective homomorphisms
$$
\mathcal O_{\bar V}(- \bar A)^{\oplus k} \subset 
\mathcal F \subset \mathcal O_{\bar V}(\bar A)^{\oplus k}
$$
for $k = \text{rank }\mathcal F$ 
which are bijective at the generic point of $\bar X$.

By Theorem B, the natural homomorphism
$\beta_*\mathcal F(m\bar L) = \pi'_*\mathcal O_{V'}(\mu^*A + m\alpha^*\bar L + mE) 
\to \pi'_*\mathcal O_{X'}((\mu^*A + m\alpha^*\bar L + mE) \vert_{X'})$ 
is surjective for any positive integer $m \ge m_2$.
Hence $X'$ is not contained in the support of $E$.

For each positive integer $m$ such that $mm_1(K_V + X)$ is
$\pi$-effective for the pair $(V, X)$, 
we construct a projective
birational morphism $\bar \mu_m: \bar V_m \to \bar V$ from a smooth variety 
such that $\bar \mu_m^*(m\bar L) = \bar Q_m + \bar N_m$ where $\bar Q_m$ is 
$\beta \circ \bar \mu_m$-free, $\bar N_m$ is effective, and that
the natural homomorphism
$(\beta \circ \bar \mu_m)_*\mathcal O_{\bar V_m}(\bar Q_m) \to 
\beta_*\mathcal O_{\bar V}(m\bar L)$ is an isomorphism.
By taking the fiber product and the normalization, we construct morphisms
$\mu'_m: V_m \to V'$ and $\alpha_m: V_m \to \bar V_m$.  
Set $\mu_m = \mu \circ \mu_m': V_m \to V$.
Let $N'_m$ be the $\pi \circ \mu_m$-fixed part of 
$\mu_m^*(mm_1(K_V + X))$.
Since $\beta_*\mathcal O_{\bar V}(m\bar L) \to \pi_*\mathcal O_V(mm_1(K_V + X))$ 
is bijective, we have $N'_m = \alpha_m^*\bar N_m + mE$.
Thus $m\bar L$ is $\beta$-effective for the pair $(\bar V, \bar X)$, 

Since the numerical Kodaira dimension of the generic
fiber of $\alpha$ is zero,
the numerical Kodaira dimension of the generic fiber of 
$\alpha \vert_{X'}$ is also zero by Theorem B  
Since $\text{dim }X_{\xi} = \text{dim }V_{\eta}$ 
and $\nu(X_{\xi}) = \nu(V_{\eta})$, 
$\bar L \vert_{\bar X}$ is $\beta$-big.
Moreover, by Theorem B again,
there exist a sufficiently large integer
$m_0$ and a global section of $\beta_*\mathcal F(m_0\bar L - 2\bar A)$ 
which induces a non-zero section over $\bar X$.
Thus there exists a global section of $\mathcal O_{\bar V}(m_0\bar L - \bar A)$ 
which does not vanish identically on $\bar X$.
Hence $\bar L$ is $\beta$-big for the pair $(\bar V, \bar X)$;
there exists an effective divisor $\bar B$ whose support does not contain 
$\bar X$ and such that $m_0\bar L \sim \bar A + \bar B$.

Let $N''_m$ be the $\pi \circ \mu_m$-fixed part of 
$\mu_m^*(mm_1(K_V + X) + A)$.
Let $\mathcal F_m = \alpha_{m*}\mathcal O_{V_m}(\mu_m^*A + m\mu^{\prime *}_mE)$.
Although $\mathcal F_m$ may be different from $\bar \mu_m^*\mathcal F/\text{torsion}$, 
the natural homomorphism $\beta_*\mathcal F(m\bar L) \to 
(\beta \circ \bar \mu_m)_*\mathcal F_m(m\bar \mu_m^*\bar L)$ is bijective.
We have
\[
\begin{split}
&\mathcal O_{\bar V_m}(\bar \mu_m^*(m\bar L - \bar A))^{\oplus k}
\subset \bar \mu_m^*\mathcal F(m\bar L)/\text{torsion} \\
&\subset \mathcal O_{\bar V_m}(\bar \mu_m^*(m\bar L + \bar A))^{\oplus k}
\subset \mathcal O_{\bar V_m}((m+m_0)\bar \mu_m^*\bar L)^{\oplus k}
\end{split}
\]
where the last inclusion is defined by $\bar B$.
Thus for $m \ge m_2$, we have
\[
\begin{split}
N''_m 
&\ge \alpha_m^*(\bar N_{m + m_0} - 
\bar \mu_m^*((m+m_0)\bar L - (m\bar L - \bar A))) + (m - m_2)E \\
&= N'_{m + m_0} - 
\bar \mu_m^{\prime *}\alpha^*(2\bar A + \bar B) - (m_0 + m_2)E \\
&\ge N'_{m + m_0} - \mu_m^*C
\end{split}
\]
for $C = (2m_0 + m_2)m_1(K_V + X)$.
Therefore, we have
$$
\mathcal I^1_{mm_1(K_V + X) + A}(- C \vert_X) 
\subset \mathcal I^1_{(m+m_0)m_1(K_V + X)} \subset \mathcal J^1_{mm_1(K_V + X)}.
$$
On the other hand,
since $K_V + X$ is $\pi$-$\Bbb Q$-effective for the pair $(V, X)$, we have
$$
\mathcal J^0_{mK_X} \subset \mathcal I^1_{m(K_V + X) + A}
$$
for any non-negative integer $m$ by Lemma 3.6. 
The rest is the same as in the proof of Theorem A.
\end{proof}

%%%%%%%%%%%%%%%%%%%%%%%%%%%%%
%%%%%%%%%%%%%%%%%%%%%%%%%%%%%
%%%%%%%%%%%%%%%%%%%%%%%%%%%%%
\section{Concluding remarks}

By Theorem 1.5 and by the base point free theorem, 
small deformations of a minimal model are always 
(not necessarily $\Bbb Q$-factorial) minimal models (cf. [KMM]).
We might ask whether a similar statement holds for global deformations.
Proposition 4.1 is on the affirmative side, but we have also a
counterexample (Example 4.2).

\vskip 1pc

{\bf Proposition 4.1.}
{\em Let $f: \mathcal X \to S$ be a proper flat algebraic morphism from a 
complex variety to a germ of a smooth curve.
Assume that the fibers $X_t = f^{-1}(t)$ have only
canonical singularities for any $t \in S$.
Let $\phi: \mathcal X \to \mathcal Z$ be a projective birational morphism over $S$
which is not an isomorphism and such that $- K_{\mathcal X}$ is $\phi$-ample.
Assume that the fibers of $\phi$ are at most $1$-dimensional.
Then the morphism restricted to the generaic fiber $\phi \vert_{X_{\eta}}$ is 
not an isomorphism.}

\vskip 1pc

\begin{proof}
Let $X_0$ be the central fiber.
For sufficiently large positive integer $m$, we have
$R^1\phi_*\mathcal O_{X_0}(mK_{X_0}) \ne 0$, while 
$R^p\phi_*\mathcal O_{X_0}(mK_{X_0}) = 0$ for $p \ge 2$.
Thus the proposition follows from the upper semi-continuity theorem combined
with Theorem 1.2'.
\end{proof}

\vskip 1pc

{\bf Example 4.2.}
Let $E = \Bbb P^d$ for an integer $d \ge 2$, and let $\mathcal X$ 
be the total space of the vector bundle $\mathcal O_E(-1)^{\oplus d}$.
Let $x_0, \ldots, x_d$ be homogeneous coordinates on $E$, and let
$\xi_1, \ldots, \xi_d$ be fiber coordinates for $\mathcal X$.
Then $t = \sum_{i=1}^d x_i\xi_i$ gives a morphism
$f: \mathcal X \to S = \Bbb C$.
The central fiber $X_0 = \{t = 0\}$ contains $E$ and 
has only one ordinary double point as singularity,
which is $\Bbb Q$-factorial if $d \ge 3$.
We have $K_{\mathcal X} \vert_E = \mathcal O_E(-1)$.
There exists a birational contraction $\phi: \mathcal X \to \mathcal Z$ 
whose exceptional locus coincides with $E$. 

\vskip 1pc

The following Example 4.3 shows that the generalizations of Theorems 1.2' 
and 1.4 to the case of varieties with log terminal singularities are false.

\vskip 1pc

{\bf Example 4.3.}
Let us consider a flip of 3-folds
with terminal singularities over a germ $\mathcal Z$:
$$
\begin{CD}
\mathcal X @>{\phi}>> \mathcal Z @<{\phi^+}<< \mathcal X^+,
\end{CD}
$$
where $- K_{\mathcal X}$ is $\phi$-ample and $K_{\mathcal X^+}$ is $\phi^+$-ample.
Let $g: \mathcal Z \to S$ be a generic projection to a germ of a smooth curve
so that the central fiber $Z_0$ is a generic hyperplane section of 
$\mathcal Z$ through the singular point.
Let $f: \mathcal X \to S$ be the induced morphism, and let 
$X_0$ be the central fiber which coincides with 
the strict transform of $Z_0$ on $\mathcal X$.

Assume that $X_0$ has only log terminal singularities.
Then so has $Z_0$ because $K_{X_0}$ is negative for $\phi$.
For example, this is the case for Francia's flip: 
(1) $\mathcal X$ has only one singularity of type $\frac 12 (1,1,1)$ and 
$\mathcal X^+$ is smooth,
(2) the exceptional loci 
$C$ and $C^+$ of $\phi$ and $\phi^+$, respectively, 
are isomorphic to $\Bbb P^1$, 
(3) the normal bundle of $C^+$ in $\mathcal X^+$ is isomorphic to
$\mathcal O_{\Bbb P^1}(-1) \oplus \mathcal O_{\Bbb P^1}(-2)$,
(4) $Z_0$ has a singularity of type $\frac 13 (1,1)$,
and (5) $X_0$ has a singularity of type $\frac 14 (1,1)$.

Since $K_{\mathcal Z}$ is not $\Bbb Q$-Cartier, $\mathcal Z$ is not log terminal.
Let $m$ be a positive integer such that $mK_{Z_0}$ is a Cartier divisor.
We take $m = 3$ for Francia's flip.
Then the natural homomorphism
$\phi_*\mathcal O_{\mathcal X}(mK_{\mathcal X}) = \mathcal O_{\mathcal Z}(mK_{\mathcal Z}) 
\to \phi_*\mathcal O_{X_0}(mK_{X_0}) = \mathcal O_{Z_0}(mK_{Z_0})$ is not surjective.
Therefore, if we compactify $\mathcal X$ suitably over $S$, 
then we obtain a counterexample
to the generalization of Theorem 1.2' for the log terminal case.

In the situation of Theorem 1.4 with $X_0$ having log terminal singularities, 
one might still ask whether the general fibers $X_t$ of $f$ have only
log terminal singularities.
This is also false.  The following example is kindly communicated 
by Professor Shihoko Ishii.
By [R1, Lemma 2.7], 
one can construct from the above example a flat deformation 
$f: \mathcal V \to B$ over a germ of a smooth curve such that 
$f^{-1}(0) \cong Z_0 \times B$ and $f^{-1}(t) \cong \mathcal Z$ 
for $t \ne 0$.
On the other hand, [Is] proved that small deformations of a log terminal 
singularity have
the same lifting property for pluricanonical forms as log terminal 
singularities although they may not be $\Bbb Q$-Gorenstein.

%%%%%%%%%%%%%%%%%%%%%%%%%%%%%
%%%%%%%%%%%%%%%%%%%%%%%%%%%%%
%%%%%%%%%%%%%%%%%%%%%%%%%%%%%
\section{Correction}

This section is a correction to this paper,
which was published in Contemporary Math. {\bf 241}(1999), 193--207.

The author would like to thank Professor Meng Chen for pointing out the mistake.
 
\vskip 1pc

Theorem A holds only under the additional condition that either 

\begin{enumerate}
\item $m > 1$, or 

\item $m \ge 1$ and $\pi(X) \ne \pi(V)$.
\end{enumerate}

These 2 cases should replace the old condition that $m \ge 1$.

Indeed the error occurred in the proof of Lemma 3.4 and the last part of the proof of Theorem A where Lemma 3.4 is used.

In the latter part of the proof of Lemma 3.4, when $D$ is only assumed to be $\pi$-$\mathbf{Q}$-effective, 
$Y_m$ should be contained in the support of some member of $\vert kQ_m \vert$ for some $k > 0$.
This is true if the image of $Y_m$ in $S$ is contained in a Cartier divisor of the image of $V$.
Thus the additional condition that $\pi(X) \ne \pi(V)$ is sufficient.

In the last part of the proof of Theorem A, if $n > 1$, then $(n-1)(K_V+X)$ is $\pi$-big, and it is OK.
If $n = 1$, then we need to assume the condition $\pi(X) \ne \pi(V)$ in order to apply the latter part of Lemma 3.4.

\vskip 1pc

We do not need to change Theorems B because $H$ is $\pi$-ample.
Theorem C and all other theorems in \S 1 are also unchanged because the condition $\pi(X) \ne \pi(V)$ is automatically satisfied.

%%%%%%%%%%%%%%%%%%%%%%%%%%%%%%%%%%%%%%%%%%
%%%%%%%%%%%%%%%%%%%%%%%%%%%%%%%%%%%%%%%%%%
%%%%%%%%%%%%%%%%%%%%%%%%%%%%%%%%%%%%%%%%%%

Graduate School of Mathematical Sciences, University of Tokyo,
Komaba, Meguro, Tokyo, 153-8914, Japan. 

kawamata@ms.u-tokyo.ac.jp

\end{document}